\documentclass[11pt]{article}

\usepackage{amsmath,amsthm,amsfonts,latexsym,amssymb}

\newcommand{\IR}{{\mathbb R}}
\newcommand{\IN}{{\mathbb N}}

\title{ Existence and bounds of positive solutions \\ for a nonlinear Schr\"{o}dinger system}
\author{Benedetta {\sc Noris}\footnote{ Partially supported by MIUR, Project ``Metodi Variazionali ed Equazioni Differenziali Non Lineari". } \ \ and \ \ Miguel\ {\sc
Ramos}}
\date{\today}

\newtheorem{thm}{Theorem}[section]
\newtheorem{lem}[thm]{Lemma}

\newtheorem{prp}[thm]{Proposition}
\newtheorem{rem}[thm]{Remark}

\begin{document}
\baselineskip0.65cm

\maketitle

\noindent{\bf Abstract.} We prove that, for any $\lambda\in \IR$, the system $
-\Delta u +\lambda u = u^3-\beta uv^2$, $ -\Delta v+\lambda v =v^3-\beta vu^2$, $ u,v\in H^1_0(\Omega),$ where $\Omega$ is a bounded smooth domain of $\IR^3$, admits a bounded family of positive solutions $(u_{\beta}, v_{\beta})$ as $\beta \to +\infty$. An upper bound on the number of nodal sets  of the weak limits of $u_{\beta}-v_{\beta}$ is also provided. Moreover, for any  sufficiently large fixed value of $\beta >0$ the system admits infinitely many positive solutions.



\setcounter{section}{0}
\section{Introduction}
\label{Introduction}

We consider systems of the form
\begin{equation}\label{system_initial}
-\Delta u +\lambda u = u^3-\beta uv^2, \qquad -\Delta v+\mu v =v^3-\beta vu^2, \qquad u,v\in H^1_0(\Omega),\end{equation}
where $\Omega$ is a smooth bounded domain in $\IR^3$ and $\lambda, \mu, \beta$ are real parameters. We are mainly interested in the case where $\lambda=\mu$ and $ \beta$ is positive and large. Such a system arises when searching for standing wave solutions of the associated time dependent Schr\"odinger system, which consists of two coupled Gross--Pitaevskii equations. This has been proposed as a mathematical model to describe both phenomena arising in nonlinear optics (see for example the references in \cite{sirakov}) and binary Bose--Einstein condensation.

In this second case  the parameter $\beta$ represents the interspecies scattering length, which determines the interaction between unlike particles: the choice $\beta>0$ corresponds to repulsive interaction (we refer to \cite{parkins_walls} for an exhaustive  physical review on B--E condensation). It has been showed experimentally that a large interspecies scattering length may induce the interesting phenomenon of phase separation, that is the two different states may repel each other and form segregated domains. Hence the analysis of system (\ref{system_initial}) with $\beta$ positive and large, besides its mathematical significance, also assumes physical relevance and has recently raised a lot of interest (see references hereafter). As it concerns the self interaction of one single state, we concentrate on the focusing case (attractive interaction), which corresponds to our choice of sign in the pure power nonlinear terms in the system \eqref{system_initial}. We stress that in our results below we could replace the constants $\lambda$ and $\mu$ by trapping potentials $\lambda(x)$ and $\mu(x)$, provided these are smooth and bounded in the $C^1$-norm in $\overline{\Omega}$.

The existence of minimal energy solutions of (\ref{system_initial}) in the whole space $\IR^3$ has been established in \cite{ambrosetti_colorado, lin_wei1, maia_montefusco_pellacci, sirakov}. These results concern the focusing case, both for the attractive and repulsive problems (see also the references therein and \cite{lin_wei2,lin_wei3} for a related problem).

Namely, for definiteness let us consider the system \eqref{system_initial}   with $\lambda=\mu=1$ and solutions $u,v\in H^1(\IR^3)$; we denote by $I_{\beta}$ the associated energy functional, whose expression is given below, and 
$$
c_{\beta}=\inf\{ I_{\beta}(u,v): u\neq 0, v\neq 0, I'_{\beta}(u,v)(u,0)=0, I'_{\beta}(u,v)(0,v)=0\}.$$
In \cite[Theorem 1]{lin_wei1} it is shown that $c_{\beta}$ is not attained in  case $\beta>0$ while in \cite[Theorem 1]{sirakov} it is proved that $c_{\beta}$ is indeed attained in case $\beta \leqslant 0$, $\beta\neq -1$; in fact, $c_{\beta}$ is attained by the diagonal pair   $(\frac{w}{\sqrt{1-\beta}}, \frac{w}{\sqrt{1-\beta}})$ where $w$ is a positive ground state solution of the equation $-\Delta w + w=w^3$ in $H^1(\IR^3)$. The existence of a positive solution was already observed in \cite[Thorem 5.4]{ambrosetti_colorado} in case $\beta\in (-1,0)$. Moreover, by combining \cite[Theorem 1]{sirakov} with  
 \cite[Theorem 5.3]{ambrosetti_colorado} or \cite[Theorem 2.3]{maia_montefusco_pellacci} we infer that 
  if  $\beta <-1$ then $c_{\beta}=I_{\beta}(\frac{w}{\sqrt{1-\beta}}, \frac{w}{\sqrt{1-\beta}})$ is the least energy critical level among all non-zero solutions of the problem, i.e. solutions $(u,v)$ with $u\neq 0$ {\em or} $v\neq 0$.  In contrast with the quoted negative result in \cite{lin_wei1}, in case $\beta>0$ it is also proved in \cite[Theorem 2]{sirakov} that a solution  with non-zero and non-negative components always exists; the existence of a non-zero solution was already observed in \cite[Theorem 5.4]{ambrosetti_colorado}. All  these solutions are shown to be radially symmetric. In connection with our Theorem  \ref{bound} below, we mention that the computation of the Morse index of the solutions is a crucial tool in the work \cite{ambrosetti_colorado}. We stress that the quoted papers deal with more general systems, namely by allowing linear terms $\lambda u$, $\mu v$ with $\lambda\neq \mu$, and nonlinear terms $\mu_1u^3$, $\mu_2v^3$ with $\mu_1\neq \mu_2$; in particular, in this case one can find ranges of  $\beta<0$ for which the problem has no solutions  with non-zero and non-negative components at all, see \cite[Theorem 1]{sirakov}.


We now concentrate on the bounded domain case. 
 We denote by $\lambda_1(\Omega)$ the first eigenvalue of  $(-\Delta,H^1_0(\Omega))$.  It has been proved by Dancer, Wei and Weth in \cite[Theorem 1.2]{dancer_wei_weth} that, for any fixed $\lambda=\mu>-\lambda_1(\Omega)$ and $\beta >0$ sufficiently large (specifically, $\beta \geqslant 1$ in the case of system \eqref{system_initial}), the system admits a positive solution $(u,v)$ (i.e. $u>0$ and $v>0$ in $\Omega$); in fact, they proved that the system  admits an unbounded sequence of positive solutions $(u_{k,\beta},v_{k,\beta})_{k\in \IN}$,  in the sense that $||u_{k,\beta}||_{H^1_0(\Omega)}+||v_{k,\beta}||_{H^1_0(\Omega)} \to \infty$ as $k\to \infty$. Their proof also provides a bound $||u_{k,\beta}||_{H^1_0(\Omega)}+||v_{k,\beta}||_{H^1_0(\Omega)} \leqslant C_k$ as $\beta\to +\infty$, for every fixed $k\in \IN$. We mention that such {\em  large} positive solutions do not exist in case $\beta < 1$ (cf. \cite[Theorem 1.1]{dancer_wei_weth}), although trivial solutions can be found by simply taking a diagonal pair $(\frac{w}{\sqrt{1-\beta}}, \frac{w}{\sqrt{1-\beta}})$ where $w$ is any positive solution of the equation $-\Delta w +\lambda w=w^3$ in $H^1_0(\Omega)$.

Since we assume $\lambda=\mu$, we rewrite our problem as
\begin{equation}\label{system}
-\Delta u +\lambda u = u^3-\beta uv^2, \qquad -\Delta v+\lambda v =v^3-\beta vu^2, \qquad u,v\in H^1_0(\Omega).\end{equation}
We recall that solutions of \eqref{system} can be seen as critical points of the $C^2$ energy functional $I_{\beta}: H^1_0(\Omega)\times H^1_0(\Omega)\to \IR$,
$$
I_{\beta}(u,v)=\frac{1}{2}\int_{\Omega}(|\nabla u|^2+|\nabla v|^2 +\lambda u^2+\lambda v^2)-\frac{1}{4}\int_{\Omega}(u^4+v^4)+\frac{\beta}{2}\int_{\Omega}u^2v^2, \quad u,v\in H^1_0(\Omega).$$
The invariance of $I_{\beta}$ with respect to the involution $(u,v)\mapsto (v,u)$ is a key ingredient in the proof of the quoted existence result in  \cite{dancer_wei_weth}.  Besides, heuristically speaking, these solutions can somehow be seen as   bifurcating from  the ground state positive solutions of the single equation $-\Delta u +\lambda u=u^3$ in $H^1_0(\Omega)$, as $\beta$ decreases from $+\infty$; we stress that such positive solutions do exist since $\lambda>-\lambda_1(\Omega)$.

Of course, the latter feature changes drastically in case $\lambda\leqslant -\lambda_1(\Omega)$. However, our previous motivation to the system \eqref{system} suggests that the mere existence of positive solutions of the system should not depend on the value   of the parameter $\lambda$. We will prove that this is indeed the case, namely that positive solutions of the system do always exist, and this will be due to the actual presence of the (sufficiently large) parameter $\beta$.

\begin{thm} \label{existence} For any $\lambda\in \IR$ and sufficiently large $\beta>0$, the system \eqref{system}
admits an unbounded sequence of solutions $(u,v)$ with $u>0$, $v>0$, $u\neq v$.
\end{thm}

In order to prove Theorem \ref{existence} we will   introduce a suitable minimax framework which takes advantage,  as in the work of \cite{dancer_wei_weth}, of the above mentioned symmetry property of the energy functional.  This is described in Section  \ref{Minimax}. We mention that in the case when the quadratic part of $I_{\beta}$  is coercive (that is, in case  $\lambda>-\lambda_1(\Omega)$), bounds on the critical points and, more generally, on the Palais-Smale sequences of $I_{\beta}$ follow immediately from any available bounds on the energy functional.  This is not the case when $\lambda\leqslant -\lambda_1(\Omega)$, and therefore the lack of compactness is an issue here.

In order to bypass this difficulty  we will work in Section \ref{Minimax} with a truncated problem. In the second part of our proof  we recover the original system by establishing {\em a priori } bounds on the solutions; since in our situation  energy estimates are useless, we  rely instead on the information of their Morse indices.  This will be presented in Section \ref{Morse}. 

It turns out that the estimates in Section \ref{Morse}  are uniform in   $\beta$, regardless of its magnitude and sign, and apply to not necessarily positive solutions of the system; see Section \ref{Morse} for the details. In particular, as a by-product of our argument we are able to derive a bound which is independent of $\beta$, as $\beta\to +\infty$.  

\begin{thm} \label{bound} For a given $k\in \IN$ (sufficiently large), as $\beta\to +\infty$ we can choose a positive solution $(u_{\beta},v_{\beta})$ of {\rm (\ref{system})} in such a way that their Morse indices (with respect to the functional $I_\beta$) are bounded by  $k$. Moreover,   $||u_{\beta}||_{H^1_0(\Omega)}+||v_{\beta}||_{H^1_0(\Omega)} \leqslant C_k$ for every  $\beta$.\end{thm}

We postpone a comment   on this result to Remark \ref{augmented_morse_index} below. We stress that the solutions obtained in Theorem \ref{bound} are genuine positive solutions of the system, in the sense that  $u_{\beta}>0$, $v_{\beta}>0$ and $u_{\beta}\neq v_{\beta}$. In particular, as already observed in \cite{dancer_wei_weth}, since $\int_{\Omega}u_{\beta}v_{\beta }(v_{\beta }^2-u_{\beta }^2)=0$,   this implies that  the components $u_{\beta}$ and $v_{\beta}$ are not ordered.

Now, by combining Theorem \ref{bound} with the Brezis-Kato estimates, one immediately deduces that the family $||u_{\beta}||_{L^{\infty}(\Omega)}+||v_{\beta}||_{L^{\infty}(\Omega)}$ is bounded uniformly in $\beta$. This allows to conclude that the solutions we found undergo the phenomenon of phase separation, which has, as we mentioned above, some physical relevance. The phase separation has been studied starting from \cite{chang_lin_lin_lin} and more recently in \cite{noris_tavares_terracini_verzini, wei_weth}; see also the pioneering papers \cite{conti_terracini_verzini1, conti_terracini_verzini2} where similar problems are analyzed. Using the results in \cite{noris_tavares_terracini_verzini} we deduce in particular that the family  $||u_{\beta}||_{C^{0,\alpha}(\overline{\Omega})}+||v_{\beta}||_{C^{0,\alpha}(\overline{\Omega})}$ is also bounded for any $\alpha\in (0,1)$ and, up to a subsequence, we have strong convergence in $H_0^1(\Omega)\cap C^{0,\alpha}(\overline{\Omega})$ to a couple of disjointly supported functions $(u,v)$ with $u\neq 0$ and $v\neq 0$. The main feature here is that the limiting domains are unknown; recent results concerning the regularity of the limiting profile and its nodal set were obtained in \cite{caffarelli_lin}. Still, in our case it remains unclear whether a limit configuration $(u,v)$ (or rather, its difference $u-v$) does satisfy a non singular differential equation.  Nevertheless, we are able to provide the following information on the limit configuration.

\begin{thm} \label{bound_on_nodal_sets} Let $(u_{\beta},v_{\beta})$ be given by Theorem {\rm \ref{bound}} and let $(u,v)$ be such that $u_{\beta}\to u $ and $v_\beta\to v$ in $H^1_0(\Omega) $ as $\beta\to +\infty$. Then the number of connected components of the set $\{u+v>0\}$ is less than or equal to $k$. \end{thm}

As mentioned before,  the proofs of Theorems \ref{existence}--\ref{bound_on_nodal_sets} will be presented in Sections \ref{Minimax} and \ref{Morse}.

\noindent{\bf Acknowledgement.} This work was completed while the first author was
visiting CMAF - Lisbon on leave  from Bicocca University. The support of both
institutions is gratefully acknowledged. The second author was supported by FCT, Funda\c c\~ao para a Ci\^encia e a  Tecnologia, Financiamento Base 2008 - ISFL/1/209.

 \section{A minimax principle}
 \label{Minimax}
 \setcounter{equation}{0}
 
 Since Theorem \ref{existence} is covered by \cite[Theorem 1.2 (a)]{dancer_wei_weth} in  case   $\lambda >-\lambda_1(\Omega)$ (see also Remark \ref{remarkdancerweiweth} hereafter), we will henceforth assume $\lambda\leqslant-\lambda_1(\Omega)$, say, in order to simplify the notations,  $\lambda=-1$. That is, we look for positive solutions of the system
 $$
-\Delta u   = u+u^3-\beta uv^2, \qquad -\Delta v  =v+v^3-\beta vu^2, \qquad u,v\in H^1_0(\Omega),$$
where $\Omega$ is such that $\lambda_1(\Omega)\leqslant 1$.  We denote $u^{\pm}:=\max\{\pm u,0\}$. Since the map $u \mapsto \int_{\Omega}(u^{+})^2$ is not  of class $C^2$  in $H^1_0(\Omega)$ and also for  later purposes of compactness (see Lemma \ref{Foreveryfixedvarepsilon,R} below), we replace the identity map $f(s)=s$ by a function which is superlinear near 0 and sublinear at infinity. For any small $\varepsilon>0$, let $f_{\varepsilon}: \IR\to \IR$ be the odd symmetric function given by $f_{\varepsilon}(s):=s^{1+\varepsilon}$ if $0\leqslant s\leqslant 1$ and $f_{\varepsilon}(s):=(1+\varepsilon)s-\varepsilon$ if $s\geqslant 1$. Then $ f_{\varepsilon} \in C^1(\IR; \IR)$ and the map $u\mapsto \int_{\Omega}F_{\varepsilon}(u^{+}) $ is $C^2 $  in $H^1_0(\Omega)$, where  $F_{\varepsilon}(s):=\int_0^{s}f_{\varepsilon}(\xi)\, d\xi$. Now,  for any large $R>0$, let  $f_{\varepsilon,R}: \IR\to \IR$ be the odd symmetric function given by $f_{\varepsilon,R}(s):=f_{\varepsilon}(s)$ if $0\leqslant s\leqslant R$ and $f_{\varepsilon,R}(s):=2\sqrt{R}\, \sqrt{s}-R$, if $s\geqslant R$. We observe that  $|f'_{\varepsilon,R}(s)|\leqslant 2$ $\forall \varepsilon, R>0$.

We look for solutions of the truncated system
 $$
-\Delta u   = f_{\varepsilon,R}(u^{+})+(u^{+})^3-\beta uv^2, \qquad -\Delta v  =f_{\varepsilon,R}(v^{+})+(v^{+})^3-\beta vu^2, \qquad u,v\in H^1_0(\Omega).$$
 Solutions of this system correspond to critical points of the $C^2$ functional $I=I_{\varepsilon, R, \beta}: H^1_0(\Omega)\times H^1_0(\Omega)\to \IR$,
$$
I(u,v)=I_{0}(u)+I_{0}(v)+\frac{\beta}{2}\int_{\Omega}u^2v^2,$$
where $I_{0}(u):=\frac{1}{2}||u||^2- \int_{\Omega}F_{\varepsilon,R}(u^{+})-\frac{1}{4}\int_{\Omega}(u^{+})^4$ and $||u||^2:=||u||^2_{H^1_0(\Omega)}=\int_{\Omega}|\nabla u|^2$.  We denote by $E_k$ the eigenspace associated to the first $k$ eigenfunctions of $(-\Delta, H^1_0(\Omega))$. In the sequel we assume that $\beta>0$ is sufficiently large ($\beta \geqslant 26$ is enough, see the proof below).

\begin{lem} \label{GivenM>0} Given $M>0$ we can find $k_0\in \IN$, independent of $\varepsilon, R$ and $\beta$, such that, for any $k\geqslant k_0$ there exists a large constant $\rho_k >0$ such that
$$
\inf\{ I(u,v): \, u-v\in E_{k-1}^{\perp},\, ||u-v||=\rho_k\} \geqslant M.$$
\end{lem}

\noindent{\sc Proof.} We recall that  $0\leqslant F_{\varepsilon,R}(s)\leqslant C_0s^2$ $\forall s\in \IR$, $\forall \varepsilon, R>0$, for some $C_0>0$.  We denote $C_1=4C_0^2|\Omega|$. Let $u-v\in E_{k-1}^{\perp}$, $ ||u-v||=\rho_k$. It is clear that, for every $k$ sufficiently large, $\rho_k$ can be chosen in such a way that, for such pairs $(u,v)$,
$$
\frac{1}{4}\int_{\Omega}|\nabla (u-v)|^2-\frac{1}{2}\int_{\Omega}(u-v)^4\geqslant M+C_1.$$
This implies 
\begin{eqnarray*}
I(u,v)&\geqslant &M+C_1+(\frac{\beta}{2}+3)\int_{\Omega}u^2v^2+\frac{1}{4}\int_{\Omega}(u^4+v^4) -C_0\int_{\Omega}(u^2+v^2)\\
& &-2\int_{\Omega}uv(u^2+v^2)\\
&\geqslant & M+C_1+\frac{\beta-26}{2}\int_{\Omega}u^2v^2+\frac{1}{8}\int_{\Omega}(u^4+v^4)-C_0\int_{\Omega}(u^2+v^2)\\
&\geqslant & M+\frac{\beta-26}{2}\int_{\Omega}u^2v^2.
\end{eqnarray*}
The conclusion follows, provided $\beta \geqslant 26$.\qed

We denote $H:=H^1_0(\Omega)\times H^1_0(\Omega)$ and by $\sigma$ the involution $\sigma(u,v)=(v,u)$. Also, for a large positive constant $R_k>\rho_k$, we let $Q_k:=B_{R_k}(0)\cap E_k$, $\partial Q_k:=\{u\in Q_k: ||u||=R_k\}$. Let
$$
\Gamma_k:=\{\gamma\in C(Q_k; H): \gamma(-u)=\sigma(\gamma(u)) \; \forall u\in Q_k \; \mbox{ and } \, \gamma(u)=(u^{+},u^{-}) \; \forall u\in \partial Q_k\}.$$
We observe that by denoting $\gamma(u)=(\gamma_1(u),\gamma_2(u))$ then $\gamma_1(-u)=\gamma_2(u) $ and $\gamma_2(-u)=\gamma_1(u)$ $\forall \gamma \in \Gamma_k$. The associated map $\theta(u):=\gamma_1(u)-\gamma_2(u) $ is therefore continuous and odd symmetric, and moreover $\theta(u)=u$ $\forall u\in \partial Q_k$. Our next lemma is then a direct consequence of the Borsuk-Ulam theorem.

\begin{lem} \label{ForeverygammanGammak} For every $\gamma\in \Gamma_k$,
$$
\gamma(Q_k)\cap S_k\neq \emptyset, \qquad \mbox{ with } S_k:=\{ (u,v)\in H: u-v\in E_{k-1}^{\perp},\, ||u-v||=\rho_k\}.$$
\end{lem} 

Now, let
$$
c_k=c_{k,\varepsilon,R,\beta}:=\inf_{\gamma\in \Gamma_k}\sup_{u\in Q_k}I(\gamma(u)).$$

By the previous two lemmas, $c_k\to+\infty$ as $k\to \infty$, uniformly in $\varepsilon$, $R$ and $\beta$. Also, it is clear that 
$$
I(u^+,u^{-})\leqslant \frac{1}{2}||u||^2 -\frac{1}{4}\int_{\Omega}u^4<0, \qquad \forall u\in \partial Q_k,$$ 
provided $R_k$ is chosen sufficiently large.

\begin{lem} \label{Foreveryfixedvarepsilon,R} For every fixed $\varepsilon$, $R$ and $\beta$, the functional $I$ satisfies the Palais-Smale condition in $H$.
\end{lem}

\noindent {\sc Proof.} Suppose $I(u_n,v_n)\leqslant C$ $\forall n$ and $I'(u_n,v_n)\to 0$ as $n\to \infty$. Then
\begin{eqnarray*}
||u_n||^2+||v_n||^2&+&\int_{\Omega}(f_{\varepsilon,R}(u_n^{+})u_n^{+}-4F_{\varepsilon,R}(u_n^{+}))+
 \int_{\Omega}(f_{\varepsilon,R}(v_n^{+})v_n^{+}-4F_{\varepsilon,R}(v_n^{+}))\\
&=&4I(u_n,v_n)-I'(u_n,v_n)(u_n,v_n)\leqslant  C+C(||u_n||+||v_n||).
\end{eqnarray*}
Since $f_{\varepsilon,R} $ is sublinear at infinity, the sequence $||u_n||^2+||v_n||^2$ is bounded. It is easy to conclude.\qed

\begin{prp} \label{Foreveryfixedvarepsilon} For every fixed $\varepsilon$, $R$, $\beta$ and $k$, with $\beta$ and $k$ sufficiently large, there exists a critical point $(u_k,v_k)$ of $I$ such that $I(u_k,v_k)=c_k$. Moreover, $(u_k,v_k)$ can be chosen in such a way that its Morse index $m(u_k,v_k)$ is less than or equal to  $k$.\end{prp}

\noindent {\sc Proof.} Since $I(\sigma(u,v))=I(u,v)$ $\forall (u,v)\in H$, the gradient flow $\eta$ associated to $\nabla I$ is $\sigma$-equi-invariant (that is, $\eta(\sigma(u,v))=\sigma(\eta(u,v))$ $ \forall (u,v)\in H$); in particular $(\eta\circ \gamma)(-u)=\sigma((\eta\circ \gamma)(u))$ $ \forall \gamma\in \Gamma_k, u\in Q_k$.  Therefore, in view of the previous lemmas (just take $M=1$ in Lemma \ref{GivenM>0}) and by using a standard argument,  $c_k$  is a critical value of the functional $I$. Since $E_k$ is a $k$-dimensional space, the statement concerning the Morse index is also classical, see e.g. \cite{ghoussoub}.\qed\\

\noindent {\bf Proof of Theorem \ref{existence}.} For the sake of clarity we split the proof into several steps.\\
{\sc Step 1.} For a given $M>0$, let us fix $k$ so large that  $c_k\geqslant M$ uniformly in $\varepsilon, R, \beta$ (this is possible by virtue of Lemma \ref{GivenM>0} and Lemma \ref{ForeverygammanGammak}) and let $(u_{\varepsilon,R,\beta},v_{\varepsilon,R,\beta})$ be given by the above proposition; we simplify the notations by dropping the subscripts,  and write $(u,v)$ in the following computations. By using the fact that $I'(u,v)(u^{-},0)=0$ and $I'(u,v)(0,v^{-})=0$ we see that $u\geqslant 0$ and $v\geqslant 0$. In this way we solve the problem
 $$
-\Delta u   = f_{\varepsilon,R}(u)+u^3-\beta uv^2, \qquad -\Delta v  =f_{\varepsilon,R}(v)+v^3-\beta vu^2, \qquad u,v\in H^1_0(\Omega), \; u,v\geqslant 0.$$
{\sc Step 2.} Since $m(u,v)\leqslant k $  there exists $C_k>0$, independent of $\varepsilon, R$ and $\beta$, such that $||u||+||v||+||u||_{L^{\infty}(\Omega)}+||v||_{L^{\infty}(\Omega)}\leqslant C_k$. We postpone to Section \ref{Morse} the proof of this fact. In particular, by choosing $R$ large enough we conclude that $(u,v)$ solves the problem
 $$
-\Delta u   = f_{\varepsilon}(u)+u^3-\beta uv^2, \qquad -\Delta v  =f_{\varepsilon}(v)+v^3-\beta vu^2, \qquad u,v\in H^1_0(\Omega), \; u,v\geqslant 0.$$
{\sc Step 3.} Since the above bound is uniform in $\varepsilon$, we can pass to the limit in the truncated system, as $\varepsilon \to 0$. This yields a limit solution, still denoted by $(u,v)$, satisfying
 $$
-\Delta u   = u+u^3-\beta uv^2, \qquad -\Delta v  =v+v^3-\beta vu^2, \qquad u,v\in H^1_0(\Omega), \; u,v\geqslant 0.$$
We stress that we have strong convergence in $H^1_0(\Omega)\times H^1_0(\Omega)$ as $\varepsilon \to 0$; in particular, the energy levels pass to the limit. Of course, the bound $C_k$ also holds in the limit, uniformly in $\beta$.\\
{\sc Step 4.} Arguing by contradiction, suppose $v\equiv 0$. Then $u$ solves $-\Delta u=u+u^3$, $u\in H^1_0(\Omega)$, $u\geqslant 0$. Since $u\neq 0$ (this is because $I(u,v)>0$), the strong maximum principle implies $u>0$ in $\Omega$. We multiply the equation by the first positive eigenfunction $\varphi_1$ of $(-\Delta, H^1_0(\Omega))$ and conclude that $\lambda_1(\Omega)\int_{\Omega}u\varphi_1=\int_{\Omega}u\varphi_1+
\int_{\Omega}u^3\varphi_1>\int_{\Omega}u\varphi_1$, hence $ \lambda_1(\Omega)>1$. This contradicts our assumption that  $\lambda_1(\Omega)\leqslant 1$. In conclusion, we have that $u,v\geqslant 0$ and $u,v\not\equiv 0$. It follows again from the strong maximum principle that $u,v>0$ in $\Omega$.\\
{\sc Step 5.} Suppose now that $v=u$. Then $u$ solves in $H^1_0(\Omega)$ the equation $-\Delta u+(\beta-1)u^3=u$ and $I(u,u)=\int_{\Omega}|\nabla u|^2-\int_{\Omega}u^2+\frac{1}{2}(\beta-1)\int_{\Omega}u^4=-\frac{\beta-1}{2}\int_{\Omega}u^4<0$. This contradicts the fact that $I(u,v)>0$. In conclusion, $v\neq u$. Since the lower bound $M$ on the energy level can be chosen  arbitrarily large, we have finished the proof of Theorem \ref{existence}.\qed

\begin{rem} \label{remark_bounded_morse_index} {\rm From the previous proof we also deduce that,  for fixed $k$, the bound on the Morse index is preserved in the limit, that is $m(u,v)\leqslant k$ for every $\beta$.  Indeed, for the moment we denote by $(u_\varepsilon,v_\varepsilon)$ the solution of the approximated problem found in Step 1 above; we stress that, as shown in   Step 2, the solution does not depend on $R$ and that the bound on the Morse index does not depend on $\beta$. We denote by  $(u,v)$ the limit of $(u_\varepsilon,v_\varepsilon)$ as $\varepsilon\to0$. Since $u$ and $v$ are positive in $\Omega$, the Lebesgue convergence theorem yields $\int_\Omega f'_\varepsilon(u_\varepsilon)\varphi^2\to \int_\Omega \varphi^2$ and the same for $v_\varepsilon$, and so
$$
I_\varepsilon''(u_\varepsilon,v_\varepsilon)(\varphi,\psi)(\varphi,\psi) \to I''(u,v)(\varphi,\psi)(\varphi,\psi) \quad \text{as }\ \varepsilon \to 0,
$$
for every $\varphi, \psi \in {\cal D}(\Omega)$, which immediately implies the claim.}
\end{rem}
\begin{rem} \label{remarkdancerweiweth} {\rm As we mentioned above, the case when $\lambda>-\lambda_1(\Omega)$ is covered by the results in \cite{dancer_wei_weth}. In this case one can use constrained minimization on the Nehari  manifold associated to the system, since the functional turns out to be coercive over this manifold. Our method provides an alternate proof of  \cite[Theorem 1.2 (a)]{dancer_wei_weth}, with the additional information on the Morse index of the solutions. We point out that in the case when $\lambda>-\lambda_1(\Omega)$  the Palais-Smale condition holds for the original functional $I$, and so there is no need  for arguing by means of a truncated problem, as we did above. As for the argument in the previous Step 4 (non-vanishing of the components of the solution pair), we can replace it by the observation that, according to a celebrated result in  \cite{gidas_spruck}, the positive solutions of the elliptic equation $-\Delta u+\lambda u =u^3$ in $H^1_0(\Omega)$ are   {\em a priori } bounded in $H^1_0(\Omega)$, whereas our solutions have arbitrarily large energy levels. }
\end{rem}

We will close the section  by establishing Theorem \ref{bound_on_nodal_sets}.

\noindent {\bf Proof of Theorem \ref{bound_on_nodal_sets}.} For a fixed and sufficiently large $k\in \IN$ and for some sequence $\beta\to +\infty$, let $(u_{\beta},v_{\beta})$ be given by Theorem \ref{bound}, so that $m(u_{\beta},v_{\beta})\leqslant k$ for every $\beta$. We will show in Section \ref{Morse} that both $(u_{\beta})$ and $(v_{\beta})$ are bounded in $H^1_0(\Omega)\cap L^{\infty}(\Omega)$ and so, as explained in the Introduction, up to a subsequence we may assume that $u_{\beta}\to u$ and $v_{\beta}\to v$ in $H^1_0(\Omega)$ as $\beta \to +\infty$. As proved in \cite{noris_tavares_terracini_verzini}, it holds that $uv=0$; moreover, $u$ satisfies the equation $-\Delta u=u+u^3$ in the open set $\{u>0\}$, and similarly for $v$. Let $\omega$ be a connected component of $\{u+v>0\}$ and denote $\overline{u}:=u|_{\omega}$,   $\overline{v}:=v|_{\omega}$; then $ \overline{u}, \overline{v} \in H^1_0(\omega)$, see e.g. \cite[Theorem IX.17 and Remark 20]{brezis}. In order to prove the theorem it would be enough to show that if $\beta$ is sufficiently large then $I''(u_{\beta},v_{\beta})(\overline{u},\overline{v}) (\overline{u},\overline{v}) <0$. We prove a slightly different version of this property which is sufficient to our purposes.

Without loss of generality, suppose $\overline{u}>0$ and $\overline{v}=0$.  Since $\overline{u}\in H^1_0(\omega)$, we can fix $\varphi\in {\cal D}(\omega)$ in such a way that
$$
\int_{\omega}|\nabla \varphi|^2-\int_{\omega}\varphi^2-3\int_{\omega}\overline{u}^2\varphi^2\leqslant \int_{\omega}|\nabla \overline{u}|^2-\int_{\omega}\overline{u}^2-2\int_{\omega}\overline{u}^4.$$
By testing the equation
$-\Delta u=u+u^3$ in $\omega$ with $\overline{u}$ we see that $\int_{\omega}|\nabla \overline{u}|^2=\int_{\omega}\overline{u}^2+\int_{\omega}\overline{u}^4$, and so
$$
\int_{\omega}|\nabla \varphi|^2-\int_{\omega}\varphi^2-3\int_{\omega}\overline{u}^2\varphi^2\leqslant -\int_{\omega}\overline{u}^4<0.$$
Now,
$$I''(u_{\beta},v_{\beta})(\varphi,0)(\varphi,0)=\int_{\omega}|\nabla \varphi|^2-\int_{\omega}\varphi^2-3\int_{\omega}u_{\beta}^2\varphi^2+\beta\int_{\omega}v_{\beta}^2\varphi^2,$$
and so, in order to prove that $I''(u_{\beta},v_{\beta})(\varphi,0)(\varphi,0)<0$ (from which the theorem follows),  it is enough to show that $\beta\int_{\omega}v_{\beta}^2\varphi^2\to 0$ as $\beta \to +\infty$.

In order to prove this, we first observe that, by using a compactness argument, it is sufficient to prove that
$\beta\int_{B_{R}}v_{\beta}^2\varphi^2\to 0$ for any ball $B_R=B_{R}(x_0)\subset B_{3R}(x_0)\subset \omega$, $x_0 \in$ supp$\varphi$. Next we observe that, over the ball $B_{2R}$, it holds that
$$
-\Delta v_{\beta}=v_{\beta}\, (1+v_{\beta}^2-\beta u_{\beta}^2)\\
\leqslant  v_{\beta}\, (C-\beta\inf_{B_{2R}} u_{\beta}^2).
$$
Since $\liminf_{\beta\to +\infty}\inf_{B_{2R}} u_{\beta}^2>0$, we deduce  that
$$
-\Delta v_{\beta} 
\leqslant  -C'\beta v_{\beta} \quad \mbox{ in }  B_{2R}, \qquad v_{\beta}\leqslant ||v_{\beta}||_{L^{\infty}(\omega)} \quad \mbox{ on }  \partial B_{2R}.$$
It follows then from \cite[Lemma 4.4.]{conti_terracini_verzini3} that $\int_{B_{R}}v_{\beta}^2\leqslant 4||v_{\beta}||^2_{L^{\infty}(\omega)} e^{-R\sqrt{C'\beta}}$ and this yields our claim.\qed

\begin{rem} \label{augmented_morse_index} {\rm The minimax class $\Gamma_k$ can be replaced by a similar one in such a way that we have the additional information $m^{*}(u_{\beta},v_{\beta})\geqslant k$, where the latter denotes the augmented Morse index of the critical point $(u_{\beta},v_{\beta})$ of the energy functional (see e.g. \cite{ghoussoub}). This fact somehow suggests that the number of connected components of the set $\{u+v>0\}$ of the limit configurations can be chosen to be arbitrarily large. This, however, remains an open problem.}
\end{rem}

\section{A priori bounds via Morse index}
\label{Morse}
\setcounter{equation}{0}

In this section we prove some estimates that were used in the proof of Theorem  \ref{existence}.  This is the content  of our next proposition (see also the subsequent remark), which, together with Remarks \ref{remark_bounded_morse_index} and \ref{remarkdancerweiweth},  also implies Theorem \ref{bound}.

In the sequel we consider a system
of the form
\begin{equation}\label{general_system}
-\Delta u +\lambda u = u^3-\beta uv^2, \qquad -\Delta v+\mu v =v^3-\beta vu^2, \qquad u,v\in H^1_0(\Omega),
\end{equation}
where $\Omega$ is a smooth bounded domain in $\IR^3$ and $\lambda, \mu, \beta$ are real parameters. Solutions of the system (not necessarily positive) are critical points of the functional $I: H^1_0(\Omega)\times H^1_0(\Omega)\to \IR$,
$$
I(u,v)=I_{\lambda}(u)+I_{\mu}(v)+\frac{\beta}{2}\int_{\Omega}u^2v^2,$$
where $I_{\lambda}(u):=\frac{1}{2}||u||^2+\frac{\lambda}{2}\int_{\Omega}u^2-\frac{1}{4}\int_{\Omega}u^4$. We observe that, for every $u,v,\varphi,\psi\in H^1_0(\Omega$),
$$
I'(u,v)(\varphi,\psi)=I_{\lambda}'(u)\varphi + I_{\mu}'(v)\psi+\beta \int_{\Omega}(uv^2\varphi+vu^2\psi)$$
and
$$
I''(u,v)(\varphi,\psi)(\varphi,\psi)=I_{\lambda}''(u)\varphi\varphi + I_{\mu}''(v)\psi\psi+\beta\int_{\Omega}(u^2\psi^2+v^2\varphi^2+4uv\varphi\psi),$$
with
$I_{\lambda}''(u)\varphi\varphi=||\varphi||^2+ \lambda\int_{\Omega}\varphi^2-3\int_{\Omega}u^2\varphi^2$.

\begin{prp} \label{Let(u_beta,v_beta} Let $(u_{\beta},v_{\beta})$ be a family  of solutions of the system {\rm \eqref{general_system}}. If the family of Morse indices $m(u_{\beta},v_{\beta})$ is bounded, as well as the coefficients $\lambda$ and $\mu$, so is the family $||u_{\beta}||+||v_{\beta}||$.
\end{prp}

\begin{rem} {\rm  It will be clear that the subsequent proof also applies for a more general system
$$
-\Delta u +f_{\lambda}(u) = u^3-\beta uv^2, \qquad -\Delta v+f_{\mu}(v) =v^3-\beta vu^2, \qquad u,v\in H^1_0(\Omega),$$
where $f_{\lambda}$, $f_{\mu}$ are $C^1$ functions in $\IR$ such that $|f'_{\lambda}(s)|+|f'_{\mu}(s)|\leqslant C_0$ $\forall s$, uniformly for bounded $\lambda, \mu$. Suppose now that  $\beta\to +\infty$. By combining Proposition \ref{Let(u_beta,v_beta}  with the Brezis-Kato estimates we deduce that  the family $||u_{\beta}||_{L^{\infty}(\Omega)}+||v_{\beta}||_{L^{\infty}(\Omega)}$ is bounded. These facts were used in the proof of Theorem \ref{existence} as presented in Section \ref{Minimax}.}
\end{rem}

\noindent {\sc Proof.} We adapt an argument  in \cite{ramos}. For simplicity of notations, we omit the subscript $\beta$ in $(u_{\beta},v_{\beta})$. We split the proof into three steps.

\noindent {\sc Step 1.} For any given vector field $V=(V_1,V_2,V_3)$, let $W$ be the Poho\u{z}aev-type vector field
$$
W=\langle \nabla u,V\rangle \nabla u-\frac{1}{2}|\nabla u|^2V+
\langle \nabla v,V\rangle \nabla v-\frac{1}{2}|\nabla v|^2V+\frac{1}{4}Q(u,v)V,$$
where
$$Q(u,v):=u^4+v^4-2\beta u^2v^2-2\lambda u^2-2\mu v^2.$$
A straightforward computation, using also the system, shows that
$$
{\rm div}W=\sum_{i,k=1}^{3}\frac{\partial u}{\partial x_i}\frac{\partial u}{\partial x_k}\frac{\partial V_i}{\partial x_k}+\sum_{i,k=1}^{3}\frac{\partial v}{\partial x_i}\frac{\partial v}{\partial x_k}\frac{\partial V_i}{\partial x_k}-\frac{1}{2}(|\nabla u|^2+|\nabla v|^2){\rm div}V+\frac{1}{4}Q(u,v){\rm div} V.$$
For a given point $x_0\in \Omega$ and smooth function $\varphi\in {\cal D}(B_r(x_0))$, with $B_r(x_0)\subset \Omega$, we let $V(x)=x$. Since $0=\int_{\Omega}{\rm div}(W\varphi^2)$ and   ${\rm div}V=3$, we deduce that
$$
\frac{3}{2}\int_{\Omega}Q(u,v)\varphi^2=\int_{\Omega}(|\nabla u|^2+|\nabla v|^2)\varphi^2 +\gamma(u,v),$$
with
$$
|\gamma(u,v)|\leqslant C\int_{\Omega}(|\nabla u|^2+|\nabla v|^2+|Q(u,v)|)|\nabla \varphi^2|.$$
A similar conclusion can be derived in case $x_0\in \partial \Omega$, provided $r$ is sufficiently small. In this case, we can choose a suitable vector field $V$ such that $|| DV-Id ||_{L^{\infty}(B_r(x_0))}={\rm o}(1)$ as $r\to 0$ and $\langle V(x),\nu_x\rangle=0$ for every $x\in B_r(x_0)\cap \partial \Omega$; here $\nu_x$ denotes the unit outward normal of $\Omega$ at the point $x$. Moreover, this vector field has the remarkable property that its divergent is {\em constant} (see Remark \ref{constant_divergent}  hereafter), namely ${\rm div}V=3$. It follows then as above that
$$
\frac{3}{2}\int_{\Omega}Q(u,v)\varphi^2\leqslant (1+{\rm o}(1))\int_{\Omega}(|\nabla u|^2+|\nabla v|^2)\varphi^2+ \gamma(u,v),$$
with ${\rm o}(1)\to 0$ as $r\to 0$.
On the other hand, using the system we see that
$$
\int_{\Omega}(|\nabla u|^2+|\nabla v|^2)\varphi^2=\int_{\Omega}Q(u,v)\varphi^2+\lambda \int_{\Omega}u^2\varphi^2+\mu\int_{\Omega}v^2\varphi^2-\int_{\Omega}\langle \nabla(\frac{u^2+v^2}{2}),\nabla \varphi^2\rangle.$$
Combining this with the previous inequality and by using a compactness argument, we conclude that it is possible to fix a small number $r>0$ and a finite number of points $x_1,\ldots,x_{\ell}\in \overline{\Omega}$ in such a way that $\overline{\Omega}\subset \cup_{i=1}^{\ell}B_r(x_i)$ and, for any smooth function $\varphi\in {\cal D}(B_{2r}(x_i))$,
$$
\int_{\Omega}(|\nabla u|^2+|\nabla v|^2)\varphi^2\leqslant C_{\lambda,\mu}\int_{\Omega}(u^2+v^2)\varphi^2+\gamma(u,v).$$
We apply the Poincar\'e inequality to the functions $u\varphi$ and $v\varphi.$ By taking a smaller $r$ if necessary, so that the $L^2$-norms are absorbed into the left hand member, this leads to the final estimate
$$
\int_{\Omega}(|\nabla u|^2+|\nabla v|^2)\varphi^2\leqslant \gamma(u,v),$$
with
$$
|\gamma(u,v)|\leqslant C\int_{\Omega}(|\nabla u|^2+|\nabla v|^2+|Q(u,v)|)|\nabla \varphi^2|+(u^2+v^2)|\nabla \varphi|^2).$$

\noindent{\sc Step 2.} Suppose first that $\beta\geqslant 0$. Since the solutions have bounded Morse indices, in each set $B_{2r}(x_i)\setminus B_r(x_i)$ we can find an annulus $A_i=\{x: a<|x-x_i|<b\}$ in such a way that
$$
I''(u,v)(u\psi,0)(u\psi,0)\geqslant 0 \qquad \mbox{ and } \qquad 
I''(u,v)(0,v\psi)(0,v\psi)\geqslant 0,$$
where $\psi\in {\cal D}(B_{2r}(x_i))$ is such that $0\leqslant \psi\leqslant 1$ and $\psi|_{A_i}=1$. We point out that $||\nabla \psi||_{L^{\infty}(\IR^3)}$ is bounded uniformly in $\beta$; to be precise, 
$||\nabla \psi||_{L^{\infty}(\IR^3)}\leqslant Cm(u_{\beta},v_{\beta})/r$, for some universal constant $C>0$.

Now, the inequality $I''(u,v)(u\psi,0)(u\psi,0)\geqslant 0$ can be written as
$$
\int_{\Omega}|\nabla u|^2\psi^2+\lambda\int_{\Omega}u^2\psi^2+2\int_{\Omega}u\psi\langle \nabla u,\nabla \psi\rangle+\beta\int_{\Omega}u^2v^2\psi^2\geqslant 3\int_{\Omega}u^4\psi^2-\int_{\Omega}u^2|\nabla \psi|^2,$$
while it follows from the equation $-\Delta u +\lambda u = u^3-\beta uv^2$ that the left hand member above equals $\int_{\Omega}u^4\psi^2$. As a consequence,
$$
 \int_{\Omega}u^4\psi^2\leqslant \frac{1}{2} \int_{\Omega}u^2|\nabla \psi|^2.$$
 By replacing $ \psi $ with $\psi^2$ and using H\"{o}lder inequality we deduce that
 $$
 \int_{\Omega}u^4\psi^4\leqslant C,$$
 for some constant $C$ independent of $\beta$. It follows once more from the equation that
 $$
 \int_{\Omega}(|\nabla u|^2+u^4+\beta u^2v^2) \psi^4\leqslant C'.$$
 We perform  a similar computation using the inequality $I''(u,v)(0,v\psi)(0,v\psi)\geqslant 0$. This yields the final conclusion that
 $$
 \int_{A_i}(|\nabla u|^2+|\nabla v|^2+u^4+v^4+\beta u^2v^2)\leqslant C'',$$
 for some constant $C''$ independent of $\beta$. We combine this estimate with the one obtained in Step 1, by choosing $\varphi\in {\cal D}(B_{2r}(x_i))$, $0\leqslant \varphi\leqslant 1$, such that $\varphi=1$ in $B_r(x_i) $ and ${\rm supp}|\nabla \varphi | \subset A_i$. This leads to the conclusion that
 $$
 \int_{B_r(x_i)}(|\nabla u|^2+|\nabla v|^2)\leqslant C''',$$
  for some constant $C'''$ independent of $\beta$, and the proposition follows.
  
 \noindent{\sc Step 3.} Suppose now that $\beta\leqslant 0$. We replace the two first inequalities in Step 2 by the single one,
  $
  I''(u,v)(u\psi,v\psi)(u\psi,v\psi)\geqslant 0.$
This can be written as
\begin{eqnarray*}
\int_{\Omega}(|\nabla u|^2+|\nabla v|^2)\psi^2&+&\int_{\Omega}(\lambda u^2+\mu v^2)\psi^2+\frac{1}{2}\int_{\Omega}\langle \nabla(u^2+v^2),\nabla \psi^2\rangle +2\beta\int_{\Omega}u^2v^2\psi^2 \\
& \geqslant & 3\int_{\Omega}(u^4+v^4)\psi^2-4\beta\int_{\Omega}u^2v^2\psi^2 -\int_{\Omega}(u^2+v^2)|\nabla \psi|^2,
\end{eqnarray*}
while it follows from the system that the left hand member above equals $\int_{\Omega}(u^4+v^4)\psi^2$. As a consequence,
$$
 \int_{\Omega}(u^4+v^4)\psi^2+|\beta| \int_{\Omega}u^2v^2\psi^2 \leqslant C \int_{\Omega}(u^2+v^2)|\nabla \psi|^2.$$
 By replacing $ \psi $ with $\psi^2$ and using once more the system  we conclude  that
 $$
 \int_{A_i}(|\nabla u|^2+|\nabla v|^2+u^4+v^4+|\beta| u^2v^2)\leqslant C'',$$
 for some constant $C''$ independent of $\beta$.  We can finish the argument as before.\qed

\begin{rem} \label{constant_divergent} {\rm The construction and properties of the vector field $V$ used in the first step of the above proof are presented in \cite[Lemma 2.1]{ramos}, except for the claim that $V $ has constant divergent. We recall here this construction. The statements are not affected by an orthogonal change of coordinates, and therefore we can assume that $x_0=0$, $ \nu_0=(0,0,1)$ and, for a sufficiently small $r>0$, $\partial \Omega\cap B_r(0)=\{(x',x_3): x_3=\theta(x')\}\cap B_r(0)$, where $x'=(x_1,x_2)\in \IR^2$ and $\theta: \IR^2\to \IR$ is a smooth map such that $\theta(0)=0$ and $\nabla \theta(0)=0$. In this case the vector field $V$ is explicitly given by
$$
V(x',x_3)=(x',x_3+\alpha(x')), \qquad \mbox{ where } \alpha(x')=\langle \nabla \theta(x'),x'\rangle -\theta(x').$$
We observe that indeed ${\rm div}V=3$.}
\end{rem}

\begin{rem} {\rm By using again  the Poho\u{z}aev-type vector field (see Step 1 of the preceding proof) with $V(x)=\nu(x) $, the unit outward normal of $\Omega$ extended in a smooth way to the whole set $\overline{\Omega}$, we deduce that also $\int_{\partial \Omega}(|\nabla u_{\beta}|^2+|\nabla v_{\beta}|^2)\leqslant C$. }
\end{rem}

{}

\noindent Benedetta Noris\\
University of Milano-Bicocca\\
Via Bicocca degli Arcimboldi 8, 20126 Milano, Italy\\
E-mail address: {\sf b.noris@campus.unimib.it}\\

\noindent  Miguel Ramos\\
University of Lisbon, CMAF - Faculty of Science\\
Av.\ Prof.\  Gama Pinto  2, 1649-003 Lisboa, Portugal\\
E-mail address: {\sf mramos@ptmat.fc.ul.pt}

\end{document}